\documentclass{article}

\usepackage{latexsym}
\usepackage{mathtools}
\usepackage{amssymb}
\usepackage{amsmath}
\usepackage[top=2cm, bottom=3cm, left=2.5cm, right=2.5cm, headheight=0.5cm,
footskip=60pt]{geometry}%
\usepackage{subfigure}
\usepackage{graphicx}
\usepackage{color}
\usepackage{epstopdf}
\usepackage{float}
\usepackage[T1]{fontenc}
\newlength{\cqfd}
\def\R{\mathbb{R}}

\def\N{\mathbb{N}}
\def\Z{\mathbb{Z}}
\def\T{\mathbb{T}}

\usepackage[english,french]{babel}
\newcommand{\into}{\int_{\Omega}}
\newcommand{\he}{h_\varepsilon}
\newcommand{\ue}{\bu_\varepsilon}

\newcommand{\bu}{\bar{u}}
{}
\newtheorem{theorem}{Theorem}[section]

\newtheorem{e-proposition}[theorem]{Proposition}

\newtheorem{e-definition}[theorem]{Definition\rm}
\newtheorem{remark}{\it Remark\/}

\newtheorem{theoreme}{Th\'eor\`eme}[section]

\newtheorem{proposition}[theoreme]{Proposition}

\def\og{\leavevmode\raise.3ex\hbox{$\scriptscriptstyle\langle\!\langle$~}}
\def\fg{\leavevmode\raise.3ex\hbox{~$\!\scriptscriptstyle\,\rangle\!\rangle$}}

\begin{document}

\selectlanguage{english}
\title{Physical Relaxation Terms for Compressible Two-Phase Systems}

\selectlanguage{english}
\author{Didier Bresch\thanks{LAMA -- UMR5127 CNRS, Bat. Le Chablais, Campus Scientifique, 73376 Le Bourget du Lac, France, Email: \texttt{Didier.Bresch@univ-smb.fr}. Research of D.B. was partially supported by the SingFlows project, grant ANR-18-CE40-0027},
Cosmin Burtea\thanks{Universit\'e de Paris, Institut de Math\'ematiques de Jussieu-Paris Rive Gauche (UMR 7586), F-75205, Paris, France, Email: \texttt{cosmin.burtea@u-paris.fr}},
Fr\'ed\'eric Lagouti\`ere\thanks{Univ. Lyon, Universit\'e Claude Bernard Lyon 1, CNRS UMR5208, Institut Camille Jordan, F.--69622 cedex, France, Email: \texttt{lagoutiere@math.univ-lyon1.fr}}}

\maketitle
\centerline{Dedicated to the memory of Andro Mikeli\'c}

\medskip

\begin{abstract}
\selectlanguage{english}
In this note, we propose the first mathematical derivation of a  macroscopic Baer-Nunziato type system for compressible two-phase flows allowing two pressure state laws depending on the different phases. By doing so, we extend the results obtained by the first author and M. Hillairet [\it Annales ENS \rm (2019)]  to cover this important physical situation. A relaxation term in the mass fraction equation  is obtained without closure assumptions contrarily to theoretical-physics literature dedicated to mixture theory, see for instance “Thermo-Fluid dynamics of Two Fluid Flows” by M. Ishii. The relaxation parameter is linked to the viscosities of the different fluids (which may be small for applications) and the relaxed quantity is linked to the laws chosen at interfaces of the two-fluid system at a mesoscale.  As this paper is intended for a large audience, we start with two formal arguments leading to the effective system. This provides formal procedures which could be useful for people working in environmental studies or industrial applications to understand how mixture models may be derived. Then we propose two mathematical proofs: One with a continuous approach (Hoff solutions for compressible NS equations with pressure depending on two transported quantities and the associated two-scale limit)  and the other with a semi-discrete approach (ODE deduced from the discretization and its continuous limit). Finally, owing to the later approach, we describe some numerical experiments by comparing mesoscopic discretization and macroscopic discretization. This  shows that  theoretical proofs may be helpful in this topic to design  appropriate numerical schemes.

\medskip

\centerline{\bf R\'esum\'e}
\smallskip
Dans cette note, nous étendons les résultats obtenus par le premier auteur et M. Hillairet 
[\it Annales ENS \rm (2019)] pour couvrir des situations physiques importantes: Plus précisément, nous proposons la première justification mathématique d'un système de type Baer-Nunziato pour écoulements biphasés compressibles permettant deux lois d'état de pression dépendant des différentes phases et un terme de relaxation qui s'obtient sans hypothèse de fermeture contrairement aux ouvrages physiques dédiés à la théorie des mélanges comme le livre «Thermo-Fluid dynamic of Two Fluid Flows »par M. Ishii. La grandeur du paramètre de relaxation est liée aux viscosités des différents fluides (qui peuvent être petites) et la quantit\'e relax\'ee est li\'ee aux lois choisies aux interfaces du système bi-fluide à l'échelle m\'esoscopique. Comme cet article est destiné à un large public, nous commençons par deux arguments formels menant au système limite: Cela fournit des procédures formelles qui pourraient être utiles pour d'autres disciplines, par exemple, dans des études environnementales ou des applications industrielles. Nous proposons ensuite deux preuves mathématiques: l'une avec une approche continue (solutions \`a la  Hoff d'\'equations de Navier-Stokes compressible avec pression dépendant de deux quantités transportées et sa limite à deux échelles) et l'autre avec une approche semi-discrète (Syst\`eme d'ODE obtenu de la discrétisation et sa limite continue). Enfin, gr\^ace \`a ces approches, nous pr\'esentons des illustrations num\'eriques comparant la discr\'etisation m\'esoscopique et la discr\'etisation macroscopique. Ceci montre que les approches th\'eoriques peuvent permettre sur ce sujet de d\'efinir des sch\'emas num\'eriques appropri\'es.
\newpage

\tableofcontents

\newpage

\vskip 1\baselineskip \noindent
{\bf Justification d'un terme de relaxation physique pour syst\`emes bi-phases. }

\end{abstract}

\section*{Version fran\c{c}aise abr\'eg\'ee}

Dans cette note, on propose de pr\'esenter l'obtention rigoureuse, en dimension un d'espace, d'un mod\`ele de m\'elange à deux phases \`a une seule vitesse avec deux lois de pression diff\'erentes suivant la phase.  Ce travail fait suite \`a un  travail du premier auteur avec M. Hillairet (voir \cite{BrHu}) sur la justification de mod\`eles de type Baer-Nunziato avec une pression commune aux deux phases. Il demande, pour une g\'en\'eralisation \`a deux pressions, de nouveaux r\'esultats d'existence de solutions \`a la Hoff sur Navier-Stokes compressible avec pression d\'ependant de deux quantit\'es satisfaisant chacune une \'equation de transport. Il est ensuite n\'ecessaire d'introduire un param\`etre d'homog\'en\'eisation et de justifier une asympotique \`a deux \'echelles vers un mod\`ele cin\'etique sous-jacent: L'obtention du mod\`ele limite bi-phases provient alors de la caract\"erisation des mesures de d\'efauts sous hypoth\`ese initiale.

Plus pr\'ecisemment, on consid\`ere deux fluides compressibles (dont on conna\^{\i}t la loi d'\'etat pour chacun) r\'egis par les \'equations de Navier-Stokes avec chacun leur viscosit\'e suppos\'ee constante. En considérant que le m\'elange de ces fluides est la limite de situations o\`u les fluides sont s\'epar\'es par des interfaces (approche multi-fluides \`a l'\'echelle mesoscopique) mais  \`a une \'echelle $\varepsilon$ de plus en plus fine, on obtient un syst\`eme v\'erifi\'e par la limite $\varepsilon \to 0$, pour lequel on a une formule pour calculer la pression du m\'elange, ainsi qu'une \'equation pour la fraction volumique de chaque constituant avec  un terme de relaxation totalement justifi\'e faisant apparaître les diff\'erentes lois de pressions. Nous nous focalisons dans cette note sur le cas \`a deux composants mais le r\'esultat se g\'en\'eralise au cas à plusieurs composant sans complication.

Pour commencer la note, nous pr\'esentons en section \ref{Section_Modelization} les mod\`eles math\'ematiques qu'ils soient m\'esoscopique ou macroscopique ainsi que les objectifs de preuves
pour lier le mod\`ele macroscopique bi-phases \eqref{TwoP} au syst\`eme mesoscopique bi-fluides \eqref{meso}--\eqref{PressionBi} avec les donn\'ees initiales \eqref{eq_initialdata_scheme}. 
 Nous pr\'esentons ensuite en section \ref{Formal_asymptotics}  deux approches formelles diff\'erentes (continue ou discr\`ete) permettant d'obtenir le mod\`ele limite. Nous commen\c cons par une approche continue au travers d'une analyse \`a deux \'echelles de type WKB. Nous continuons ensuite par une approche discr\`ete qui serait plus proche de ce que nous pouvons rencontrer dans le cadre d'applications. Ces r\'esultats pr\'ecisent les calculs formels que l'on peut trouver par exemple dans \cite{SaAb}, \cite{MuGu}, \cite{AbZa} en \'evitant l'hypoth\`ese de fermeture formelle pour le terme de relaxation. Cette partie peut avoir un r\'eel int\'er\^et pour le lecteur int\'eress\'e par une compr\'ehension formelle de l'obtention du mod\`ele de m\'elange \`a partir du mod\`ele m\'esoscopique avec interfaces.
 
Nous proposons ensuite dans la section \ref{Mathematical_results} deux approches th\'eoriques de justification math\'ematiques totalement compl\'ementaires: I) Une approche d'homog\'en\'eisation sur le mod\`ele mesoscopique sous sa forme continue o\`u un th\'eor\`eme de type existence de solutions \`a la Hoff (pour NS compressible avec pression d\'ependant de deux paramètres transportés) et contr\^ole uniforme en le petit param\`etre  d'homog\'en\'eisation est n\'ecessaire pour le passage \`a la limite. II) Une approche d'homog\'en\'eisation sous un angle plutôt discret qui consiste \`a consid\'erer le syst\`eme semi-discr\'etis\'e associ\'e, d'en d\'eduire un th\'eor\`eme d'existence locale pour le syst\`eme d'ODE, de montrer un contrôle uniforme de versions discr\`etes des contr\^oles \`a la Hoff puis de passer \`a la limite vers le mod\`ele continu justifiant ainsi une analyse formelle par un des auteurs dans \cite{La}. Cette partie peut int\'eresser le lecteur motiv\'e par une justification rigoureuse au travers d'une approche continue ou d'une approche discr\'etis\'ee.

 Nous terminons  ensuite par la section \ref{Numerical_illustrations} o\`u nous montrons comment,  gr\^ace aux r\'esultats th\'eoriques, nous pouvons d\'efinir de bons sch\'emas num\'eriques et nous illustrons le tout par des résultats de simulation.  
 Cette note  est \'ecrite pour pr\'esenter aux lecteurs, qu'ils soient math\'ematiciens, physiciens ou ing\'enieurs un nouveau cadre permettant d'aborder formellement, math\'ematiquement et num\'eriquement l'obtention de certains mod\`eles multiphases \`a partir de mod\`eles multi-fluides mesoscopiques. 
 Ce travail doit \^etre vu comme une partie de la monographie compl\`ete \cite{BrBuHiLa}  sur les syst\`emes multiphasiques avec deux pressions   que les auteurs avec M. Hillairet se proposent de r\'ediger. Cette monographie rappellera les travaux pr\'ec\'edents th\'eoriques (de  \cite{BrHi2}, \cite{BrDeGhGrHi}, \cite{BrHu}, \cite{BrHu}, \cite{AmZl}, \cite{PlSo},  \cite{Se},  \cite{E}  par exemple) ainsi que les travaux de mod\'elisation et aspects num\'eriques importants  (\cite{SaAb}, \cite{DeGa}, \cite{MuGu}, \cite{AbZa}, \cite{Co} et r\'ef\'erences contenues par exemple). Elle donnera \'egalement le d\'etail des preuves et discutera quelques extensions possibles.  

\selectlanguage{english}

\section{Introduction.}
\label{sec2}
\noindent In this note, we propose to present the rigorous justification, in one space dimension, of a single velocity two-phase mixing model with two different pressure laws depending on each phase. This work follows a work by the first author with M. Hillairet (see \cite{BrHu}) on the justification of Baer-Nunziato type models with a pressure common to both phases. This asks, for a generalization to two possible different pressure state laws, new results of existence of solutions \`a la Hoff on compressible Navier-Stokes with pressure depending on two densities each satisfying a transport equation. It is then necessary to introduce a homogenization parameter and to justify a two scales asymptotic towards an underlying kinetic model: The obtaining of the two-phase limit model then comes from the characterization of the defect measures under initial hypothesis.
   More precisely, we consider two compressible fluids (of which we know the pressure state law for each) governed by the Navier--Stokes equations, each with their constant viscosity that we can write as a single compressible system with a pressure dependent on the averaged-density and a color function that is used to distinguish the fluid phase, both of which are advected by the flow.
   By considering that the mixture of these fluids is the limit of situations where the fluids are separated by interfaces (multi-fluid approach at the mesoscopic scale) but at an $\varepsilon$ scale more and more fine, we obtain a system verified by the limit $\varepsilon \to 0$, for which we have a formula to calculate the pressure of the mixture, as well as an equation for the volume fraction of each component with a completely justified relaxation term showing the different laws of pressures. We focus in this note on the two-component case but the result generalizes to the multi-component case.
   
   To start the note, we present in section \ref{Section_Modelization} the mathematical models whether they are mesoscopic or macroscopic as well as the proof objectives to link the two-phase macroscopic model \eqref{TwoP} to the bi-fluid mesoscopic system \eqref{meso} -- \eqref{PressionBi} with the initial data \eqref{eq_initialdata_scheme}.
 We then present in section \ref{Formal_asymptotics} two different formal approaches (continuous or discrete) allowing to obtain the limit model. We start with a continuous approach through a two-scale analysis of WKB type. We then continue with a discrete approach which would be closer to what we may encounter in applications. These results render mathematically rigorous the formal computations that can be found for example in \cite{SaAb}, \cite{MuGu}, \cite{AbZa} by avoiding the formal closure hypothesis for the relaxation term. This part may have a real interest for the reader interested in a formal understanding of obtaining the mixture model from the mesoscopic model with interfaces.
 
We then propose in the section \ref{Mathematical_results} two completely complementary theoretical approaches to mathematical justification: I) A homogenization approach on the mesoscopic model under its continuous form where a theorem of type existence of solutions \`a la Hoff (for  compressible Navier-Stokes with pressure depending on two transported parameters) and uniform control in the small parameter  is necessary for the passage to the limit. II) A homogenization approach from a rather numerical analysis-angle which consists in considering the associated semi-discrete system, hence to deduce a local existence theory for the ODE system, to show a uniform control of discrete versions of the energy estimates à la Hoff then to pass to the limit towards the continuous model thus justifying a formal analysis by one of the authors in \cite{La}. This part is oriented towards the reader interested in a rigorous justification through a continuous approach or a discrete approach.

Owing to the theoretical results, in the last Section, we show how to define appropriate numerical schemes and we present the results of some numerical experiments which better illustrate our results.
 This note is addressed to a wide audience and aims at presenting a new framework allowing to approach formally, mathematically and numerically the derivation of some multiphase models with different pressure laws from mesoscopic descriptions. 
 
 This work should be seen as an introduction to the complete monograph \cite{BrBuHiLa} that the authors propose to write in collaboration with M. {\sc Hillairet}. This monograph will recall previous important theoretical results (for instance  \cite{AmZl},  \cite{Se}, \cite{E}, \cite{PlSo}, \cite{BrHi2}, \cite{BrDeGhGrHi}, \cite{BrHu}, \cite{BrHu} for example) as well as results dealing with mathematical modelling and important numerical aspects (\cite{SaAb}, \cite{DeGa}, \cite{ MuGu}, \cite{AbZa}, \cite{Co} and references therein for example). It will also give full details of the proofs and discuss some possible extensions.

\section{Mathematical models and objectives.}\label{Section_Modelization}
\subsection{The mesoscale system}

In this part, we detail the process towards the mathematical justification of
the Baer-Nunziato model for multiphase mixture with a physical relaxation
term.  

 Our main assumption is to consider a {\em two-fluid mixture}
(multiphase fluid) as the {\em limit of a sequence of situations where the two fluids are separated
by sharps interfaces} (this configuration is sometimes referred as a
multifluid material). At a physical level this translates as follows: we assume that it is possible to zoom-in in the mixture at a very fine scale, which is referred to as being the mesoscopic scale.{\em We assume that at this scale the two constituents are separated which allows us to model their behaviour with the well-known Navier-Stokes system plus an equation describing the evolution of the interface}. In order to recover a macroscopic effective model, we have to zoom-out which mathematically is translated as a propagation of oscillations problem and its averaging.

Then, our second assumption is that the multifluid is governed by a
Navier-Stokes type system where the viscosity and the pressure depend on
the fluid, dependence that can be modelled by introducing a binary color parameter, which takes its values depending on the phase. More precisely,
we consider a mixture of two compressible fluids, which will be referred  in
the following as $+$ and $-$ occupying a domain $\Omega$. In the
following, physical quantities that characterize the phase $+/-$ will be
explicitly denoted by a $+/-$ lower script.

Let us denote by $\Omega_{+}\left(  t\right)$ the volume occupied by fluid $+$
and by $\Omega_{-}\left(  t\right)$ the volume occupied by fluid $-$ at
time $t$. We assume that initially we have
$\Omega_{+}\left(  0\right) \cap \Omega_{-}\left(0\right) =
\emptyset$ and $\Omega_{+}\left(0\right)
\cup \Omega_{-}\left(  0\right) = \Omega$. Lest us denote by $c(t,\cdot)$ the characteristic function of $\Omega_{+}\left(0\right)$ (and note that thus $c$ is also the volume fraction or the mass fraction of fluid $+$).
Let us denote $\rho_{+}$ and $\rho_{-}$ the densities of the two phases. The
densities are initially defined in $\Omega_{+}\left(  t\right)  $ respectively
in $\Omega_{-}\left(  t\right)$. We extend $\rho_{+}$ and $\rho_{-}$ by $0$
in $\Omega$ and slightly abusing the notation we will still call them
$\rho_{\pm}$. We see that owing to the separation hypothesis and the
definition of $c$ we have that
\begin{equation}
\rho = c\rho_{+}+\left(  1-c\right)  \rho_{-}
\label{rho=rho_++rho_-}%
\end{equation}
and
\begin{equation}
\rho_{+}=c\rho\text{, }\rho_{-}=\left(  1-c\right)  \rho. \label{rho_pm}%
\end{equation}
a.e. in $\Omega$.

We assume that fluid $+$ has a constant viscosity $\mu_+ > 0$ and
fluid $-$ has a constant viscosity $\mu_- > 0$, and denote by
$\mu(c)$ the viscosity of the multifluid:
\[
\mu(c)=c \mu_{+}+(1-c)\mu_{-}.
\]
We assuming that both fluids are barotropic fluids with
pressure laws $p_+(\rho)$ and $p_-(\rho)$, and denote by $p(c,\rho)$
the pressure law of the multifluid:
\[
p(c,\rho)=c\,p_{+}(\rho)+(1-c)p_{-}(\rho).
\]
We denote by $u(t,x)$ the velocity at time $t$ and position $x$.
The above notations and definitions lead to model the multifluid Cauchy problem
with the following mesoscopic system:

\begin{equation}
\left\{
\begin{array}
[c]{l}%
\partial_{t}c+u\,\partial_{x}c=0\hbox{ with }c\,(1-c)=0,\\
\partial_{t}\rho+\partial_{x}(\rho u)=0,\\
\partial_{t}(\rho u)+\partial_{x}(\rho u^{2})-\partial_{x}(\mu(c)\partial
_{x}u)+\partial_{x}p(c,\rho)=0,\\
\mu(c)=(c\,\mu_{+}+(1-c)\mu_{-}),\qquad p(c,\rho)=c\,p_{+}(\rho)+(1-c)p_{-}%
(\rho).\\
\rho_{|t=0}=\rho_{0},\qquad c_{|t=0}=c_{0} \hbox{ with } c_0 (1-c_0)=0,\qquad u_{|t=0}=u_{0},
\end{array}
\right.  \label{System}%
\end{equation}

Let take a moment to resume what has been done up to the present. We made a
number of hypothesis that allowed us to propose a system of equations
governing the evolution of a mixture at a mesoscopic scale. Although very
interesting at a mathematical level this system is very unlikely to be of use
in practical applications. Thus, the necessity to un-zoom back at macrospic
scale comes naturally. Using similar systems as $\left(  \text{\ref{System}%
}\right)  $, performing a formal procedure of averaging and assuming some
closure assumptions, several calculations may be encountered in
physics-papers such as \cite{DrPa}, \cite{Is} or modelling or numerical papers
such as \cite{SaAb}, \cite{AbZa}. We describe in the following
lines how to obtain mathematically equations at a macroscopic scale. The key
observation from \cite{BrHi1,BrHi2} is that when un-zooming, an observer will witness
rapid oscillations between the zones occupied by the two phases. Thus, a
multi-phase fluid can be represented by the limit of solutions that widely
oscillate in space. More precisely, following the formalism introduced in
\cite{BrHu}, \cite{BrHi1}, we obtain a macroscopic effective system for multiphase
fluids by introducing a parameter $\varepsilon$ describing the oscillation
scale. The idea is to consider a sequence of solutions of $\left(
\text{\ref{System}}\right)  $ generated by a sequence of initial data widely
oscillating in space and to analyse its limiting behaviour. A typical example of such initial data is given by
\begin{equation}
c_{0}^{\varepsilon}(x)=c_{0}(x/\varepsilon),\qquad\rho_{0}^{\varepsilon
}(x)=c_{0}(x/\varepsilon)\rho_{+,0}(x)+(1-c_{0}(x/\varepsilon))\rho
_{-,0}(x)\label{eq_initialdata_scheme}%
\end{equation}
where $c_{0}:\Omega\rightarrow\{0,1\}$ is a fixed profile and
$\rho_{0,+}$, $\rho_{0,-}$ are bounded initial data and where the initial
velocity field $u_{0}^{\varepsilon}(x)=u_{0}(x)\in H^{1}(\Omega)$. Thus,
one considers $(c^{\varepsilon},\rho^{\varepsilon},u^{\varepsilon
})_{\varepsilon>0}$ a sequence of solutions (in the sense of Hoff) of the bifluid system mentioned previously i.e.
\begin{equation}
\left\{
\begin{array}
[c]{l}%
\partial_{t}c^{\varepsilon}+u^{\varepsilon}\partial_{x}c^{\varepsilon}=0,\\
\partial_{t}\rho^{\varepsilon}+\partial_{x}(\rho^{\varepsilon}u^{\varepsilon
})=0,\\
\partial_{t}(\rho^{\varepsilon}u^{\varepsilon})+\partial_{x}(\rho
^{\varepsilon}{u^{\varepsilon}}^{2})-\partial_{x}(\mu(c^{\varepsilon}%
)\partial_{x}u^{\varepsilon})+\partial_{x}p(c^{\varepsilon},\rho^{\varepsilon
})=0,\\
\mu(c^{\varepsilon})=(c^{\varepsilon}\mu_{+}+(1-c^{\varepsilon})\mu_{-}),\\
p(c^{\varepsilon},\rho^{\varepsilon})=c^{\varepsilon}p_{+}(\rho^{\varepsilon
})+(1-c^{\varepsilon})p_{-}(\rho^{\varepsilon}),
\end{array}
\right.  \label{meso}%
\end{equation}
with
\begin{equation}
\left\{
\begin{array}
[c]{l}%
c^{\varepsilon}|_{t=0}=c_{0}^{\varepsilon}\mbox{ such that }c_{0}%
^{\varepsilon}(1-c_{0}^{\varepsilon})=0,\\
\rho^{\varepsilon}|_{t=0}=\rho_{0}^{\varepsilon},\\
u^{\varepsilon}|_{t=0}=u_{0}^{\varepsilon},
\end{array}
\right.  \label{mesoCauchy}%
\end{equation}
where the initial data are supposed to oscillate wildly in space.
The property $c_{0}^{\varepsilon}(1-c_{0}^{\varepsilon})=0$ is important
because it implies that $c^{\varepsilon}(1-c^{\varepsilon})=0$ at any time and
that thus guarantees that it is legal to compute the pressure as
\begin{equation}\label{PressionBi}
p(c^{\varepsilon},\rho^{\varepsilon})=c^{\varepsilon}p_{+}(\rho^{\varepsilon
})+(1-c^{\varepsilon})p_{-}(\rho^{\varepsilon}),
\end{equation}
the phases are \textquotedblleft pure\textquotedblright\ at any point of
$\Omega^{\varepsilon}$. The macroscopic bi-phase model is then derived letting
$\varepsilon$ tend to $0$ and computing the limit system. We prove that, if
the Cauchy data \eqref{mesoCauchy} converges weakly to $(\alpha_{0},\rho
_{0},u_{0})$ (note that the property $\alpha_{0}(1-\alpha_{0})$ is of course
lost), then, up to subsequence $(c^{\varepsilon},\rho^{\varepsilon
},u^{\varepsilon})$ converges to $(\alpha,\rho,u)$ such that
\begin{equation}
\left\{
\begin{array}
[c]{l}%
\partial_{t}\alpha+u\,\partial_{x}\alpha=\dfrac{\alpha\left(  1-\alpha\right)
}{\alpha\mu_{-}+(1-\alpha)\mu_{+}}(F_{\pm}-F_{\mp}),\\
\partial_{t}\rho+\partial_{x}(\rho u)=0,\\
\partial_{t}(\rho u)+\partial_{x}(\rho u^{2})-\partial_{x}(\mu_{\mathrm{eff}%
}\partial_{x}u)+\partial_{x}p_{eff}=0,\\
\rho=\alpha\rho_{+}+(1-\alpha)\rho_{-},\\
\mu_{\mathrm{eff}}=\dfrac{\mu_{+}\mu_{-}}{\alpha\mu_{-}+(1-\alpha)\mu_{+}},\\
p_{\mathrm{eff}}=\dfrac{\alpha p_{+}(\rho_{+})\mu_{+}+(1-\alpha)p_{-}(\rho
_{-})\mu_{-}}{\alpha\mu_{-}+(1-\alpha)\mu_{+}},\\
F_{\pm}=-\mu_{\pm}\partial_{x}u+p_{\pm}(\rho_{\pm}).
\end{array}
\right.  \label{TwoP}%
\end{equation}
Observe that the resulting homogeneized system contains only physical
quantities such has the two viscosities, the two pressure laws and the volume
fraction. We did not suppose a closure assumption with a relaxing parameter as
usually in the physical literature.

In the sequel, we present two approaches in order to derive solutions of the
bifluid system and to obtain at the limit the same macroscopic equations namely

\begin{itemize}
\item A so called continuous approach which follows the ideas introduced by
D. Bresch and M. Hillairet in \cite{BrHu}, \cite{BrHi1}, \cite{BrHi2}.
This consist in using the techniques introduced by D. Hoff in order to
construct solutions of the bifluid system. The main novelty here is the
presence of the new unknow $c$ allowing to take into account different
pressure laws. The bifluid system is a compressible Navier-Stokes system with
a viscosity given in terms of $c$ and a pressure law depending on $c$ and
$\rho$.

\item A semi-discretized approach which will render rigorous the result by one
of the authors (see \cite{La}). This approach is interesting because it allows
to numerically compute solutions for the macrospic mixture model $\left(
\text{\ref{TwoP}}\right)  $ with a numerical scheme designed on the mesoscopic
system $\left(  \text{\ref{meso}}\right)  $.
\end{itemize}

\bigskip

In Section \ref{Formal_asymptotics}, we show that we may obtain the system $\left(
\text{\ref{TwoP}}\right)  $ from the system $\left(  \text{\ref{meso}}\right)
$ first by a formal two-scale analysis at a continuous level and then by a scale analysis at a discrete level . In Section \ref{Mathematical_results}, we present the mathematical results justifying the derivation from the bifluid description to the biphase system: continuous approach and then semi-discrete approach.  In the last section \ref{Numerical_illustrations}, we illustrate the results stated in the paper. More precisely, we design two numerical schemes: one to approximate the mesoscopic system, that is to say System \eqref{meso} with a Cauchy datum \eqref{mesoCauchy}, and one to approximate the macroscopic system \eqref{TwoP} with any Cauchy datum.

\section{Formal asymptotics}\label{Formal_asymptotics}

\subsection{How to derive the relaxed equation \eqref{TwoP}${}_1$ without closure assumptions?}

\subsubsection{Continuous level -- Two scale analysis}
In the homogenization process, two-scale asymptotic is a natural tool (see for instance \cite{Al}, \cite{Se} and \cite{E}
and references cited therein). Let us recall formally, how it may be used in the compressible setting
using some ellipticity properties of the effective flux and using the renormalized approach for the transport equation. Such formal computation has been performed in \cite{BrHi2} but for the reader's convenience we rewrite here such calculation. In the bifluid setting, to get the asympotic system we formally assume the following ansatz
$$ c\, (t,x)= c \, (t,\frac{t}{\varepsilon}, x, \frac{x}{\varepsilon})$$
$$\rho(t,x) = c\, (t,\frac{t}{\varepsilon}, x, \frac{x}{\varepsilon}) \rho_+^\varepsilon(t,x)
     + \bigl(1- c \,(t,\frac{t}{\varepsilon}, x, \frac{x}{\varepsilon})\bigr) \rho_-^\varepsilon(t,x)$$
$$u(t,x) = u_0(t,\frac{t}{\varepsilon}, x, \frac{x}{\varepsilon})
         + \varepsilon u_1(t,\frac{t}{\varepsilon}, x, \frac{x}{\varepsilon})
         + \varepsilon^2 u_2(t,\frac{t}{\varepsilon}, x, \frac{x}{\varepsilon})
         + O(\varepsilon^3)
$$
assuming
$$\rho_\pm^\varepsilon (t,x) = \rho^0_\pm (t,x) + O(\varepsilon), \qquad
   c(t,\tau,x,y) \in \{0,1\} \hbox{ \it a.e. }.$$
   Plugging these informations in the bifluid system, we first get
 \begin{equation} \label{First} \partial_\tau c + u_0 \partial_y c = 0
 \end{equation}
 and
 \begin{equation}\label{Second}
 \partial_t c  + u_0 \partial_x c + u_1 \partial_y c = 0.
 \end{equation}
 The first equation provides the behavior of $c$ on a cell. This equation is
 compatible with the assumption that $c$ is an indicator function. Averaging with
respect to the fast variable $y$ the second equation, we get the following equation
\begin{equation} \label{Third}
\partial_t \alpha + \overline{ u_0\partial_x c} = - \overline{u_1 \partial_y c}
\end{equation}
with $\alpha = \overline c$ where we denote temporarily with the bar average with respect
to $y$ on a cell. Let us now quickly recall the different steps. Plugging now the ansatz in
the momentum equation we get
$$\partial_y ((c \mu_+ + (1-c) \mu_-)\partial_y u_0) = 0.$$
Multiplying this equation and integrate with respect to the $(t,\tau,x,y)$ and using that
$$ c\mu_+ + (1-c)\mu_- \ge \min (\mu_+,\mu_-)> 0$$
we get that
$$\partial_y u_0 = 0$$
and therefore coming back to the equation on $c$ and averaging with respect to $y$ that $\alpha$ does not depend on $\tau$. Looking now the main part of the momentum equation, we get
\begin{multline*}
\rho^0 (\partial_\tau u_0 + u_0 \partial_y u_0) + \partial_y p^0
    =  \partial_y((c \mu_+ + (1-c) \mu_-) \partial_y u_1) \\
       + \partial_x((c \mu_+ + (1-c) \mu_-) \partial_y u_0)
       + \partial_y((c\mu_+ + (1-c) \mu_-) \partial_x u_0).
\end{multline*}
Multiplying by $\partial_\tau u_0$ and integrating with respect to the $(t,\tau,x,y)$,
 we get that $\partial_\tau u_0 = 0.$  Thus we get
$$\partial_y((c \mu_+ + (1-c) \mu_-)\partial_y u_1) - \partial_y p^0
      = -\partial_y((c\mu_+ + (1-c) \mu_-) \partial_x u_0).$$
Denoting $\mu = c \mu_+ + (1-c) \mu^-$, this gives
$$  \mu \partial_y u_1 - p^0 =  \overline{\mu \partial_y u_1}
      - \overline{p^0}  + (\overline{\mu}- \mu) \partial_x u^0.$$
Note that using the expression of $\mu$ and the equation of $\alpha$,
we get
$$\mu \partial_y u_1 = (p^0-\overline{p^0}) - (\mu - \overline{\mu}) \partial_x u_0
                                          + (\mu_+-\mu_-) (\partial_t \alpha + u_0\partial_x \alpha)$$
which may be rewritten as
$$c\,\partial_y u_1 = \frac{c}{\mu}\Bigl((p^0-\overline{p^0}) - (\mu - \overline{\mu}) \partial_x u_0
                                          + (\mu_+-\mu_-) (\partial_t \alpha + u_0\partial_x \alpha)\Bigr)
$$
and therefore after calculation
\begin{multline*}
\overline{ c\partial_y u^1} = \frac{\alpha (1-\alpha)}{\mu_+} \bigl((p_+(\rho^0_+) - p_-(\rho_-^0)) 
       - \partial_x u_0 (\mu_+-\mu_-) + \alpha (1- \frac{\mu_-}{\mu_+}) (\partial_t + u_0 \partial_x) \alpha
\end{multline*}
using that
$$ \overline{u^1 \partial_y c} = - \overline{ c\partial_y u^1}
$$
and inserting in \eqref{Third}, we get after simple calculations
$$\partial_t \alpha + u_0 \partial_x \alpha
    = \frac{\alpha (1-\alpha )}{\alpha \mu_- + (1-\alpha) \mu_+}
    \Bigl((p_+(\rho_+^0) - p_-(\rho_-^0))  - \partial_x u_0 (\mu_+- \mu_-)\Bigr)$$

\subsubsection{Discrete level -- scale analysis}

In this section we propose a formal procedure to derive the mixture model and
we explain with very simple arguments why the volume fraction should satisfy
the equation
\[
\partial_{t}\alpha+u\partial_{x}\alpha=\frac{\alpha(1-\alpha)}{\mu}%
(p_{+}-p_{-}).
\]
We denote by $D_{t}\alpha$ the Lagrangian time derivative of $\alpha$:
\[
D_{t}\alpha=\partial_{t}\alpha+u\partial_{x}\alpha
\]
(this derivative will also be denoted $\dot{\alpha_{+}}$ in the rest of the
paper). Consider a situation where the fluids are separated (say, at a small
scale $\varepsilon$), and a point $x(t)\in{\mathbb{T}}$ at an interface
between fluid $+$ on its right and fluid $-$ on its left, for any time $t$.
Denote by $x_{+}(t)$ the center of the zone of pure fluid $+$ on the right of
$x(t)$, by $x_{-}(t)$ the center of the zone of pure fluid $-$ on the left of
$x(t)$, and
\[
\varepsilon_{+}(t)=x_{+}(t)-x(t),\varepsilon_{-}(t)=x(t)-x_{-}(t)
\]
which are supposed to be small.

We define $\alpha(t)$ by
\[
\alpha(t)=\varepsilon_{+}(t)/(\varepsilon_{+}(t)+\varepsilon_{-}(t)).
\]
Indeed this quantity represents the local (at point $x(t)$) volume fraction of
fluid $+$. Obviously one has
\[
D_{t}\varepsilon_{+}(t)=u(t,x_{+}(t))-u(t,x(t))
\]
and
\[
D_{t}(\varepsilon_{+}+\varepsilon_{-})(t)=u(t,x_{+}(t))-u(t,x_{-}(t)).
\]
This allows to write
\begin{multline}
D_{t}\alpha(t)=\frac{(\varepsilon_{+}+\varepsilon_{-})D_{t}\varepsilon
_{+}-\varepsilon_{+}D_{t}(\varepsilon_{+}+\varepsilon_{-})}{(\varepsilon
_{+}+\varepsilon_{-})^{2}}\label{dta}\\
=\frac{\varepsilon_{-}(u(t,x_{+}(t))-u(t,x(t)))-\varepsilon_{+}%
(u(t,x(t))-u(t,x_{-}(t)))}{(\varepsilon_{+}+\varepsilon_{-})^{2}}%
\end{multline}
The regularity of the solution is expected to be the following: at any time
$t$, the pressure and the space derivative of the velocity should be
continuous in space in each pure region (namely, in $(x_{-}-\varepsilon
_{-},x_{-}+\varepsilon_{-})$ and in $(x_{+}-\varepsilon_{+},x_{+}%
+\varepsilon_{+})$), but not at the point $x(t)$. At this point, what is
expected is that the effective flux $p-\mu\partial_{x}u$ is continuous (and
this continuity in space stands for the law of reciprocal forces of Newton).
In the case where the two viscosity coefficients are equal, the formal
computation is straightforward. Thus we propose to begin by assuming this
equality, and to obtain the general law for $\alpha$ in a second stage.

\begin{itemize}
\item Case where $\mu_{+} = \mu_{-} = \mu$ \newline The continuity of the
effective flux together with the regularity on pure zones expresses as
\begin{multline*}
p_{-}(t) - \mu\frac{u(t,x(t)) - u(t,x_{-}(t))}{\varepsilon_{-}} = p_{+}(t) -
\mu\frac{u(t,x_{+}(t)) - u(t,x(t))}{\varepsilon_{+}}\\
+ r(\varepsilon_{-} + \varepsilon_{+}),
\end{multline*}
where $p_{\pm}(t)$ denotes $p_{\pm}(\rho(t,x_{\pm}(t)))$ and $r$ is a function
such that $r(x) \rightarrow0$ as $x \rightarrow0^{+}$. This rewrites
\begin{multline*}
p_{+}(t) - p_{-}(t) = \mu\frac{\varepsilon_{-}(u(t,x_{+}(t)) - u(t,x(t))) -
\varepsilon_{+}(u(t,x(t)) - u(t,x_{-}(t)))}{\varepsilon_{+}\varepsilon_{-}}\\
+ r(\varepsilon_{-} + \varepsilon_{+}),
\end{multline*}
and, thanks to \eqref{dta} and letting $\varepsilon_{\pm}$ go to $0$,
\[
p_{+}(t) - p_{-}(t) = \mu\frac{(\varepsilon_{+} + \varepsilon_{-})^{2}%
}{\varepsilon_{+} \varepsilon_{-}} D_{t} \alpha_{+} = \frac{\mu}{\alpha_{+}(1
- \alpha_{+})} D_{t} \alpha_{+},
\]
which is exactly what is stated in this paper.

\item Case where $\mu_{+}\neq\mu_{-}$ \newline In this general case, it is
convenient to define the approximate space derivatives of the velocity
$d_{-}(t)$ and $d_{+}(t)$
\[
d_{-}(t)=\frac{u(t,x(t))-u(t,x_{-}(t))}{\varepsilon_{-}(t)},\quad
d_{+}(t)=\frac{u(t,x_{+}(t))-u(t,x(t))}{\varepsilon_{+}(t)}.
\]
Equipped with this, we can rewrite \eqref{dta} as
\[
D_{t}\alpha(t)=\frac{\varepsilon_{-}\varepsilon_{+}}{(\varepsilon
_{-}+\varepsilon_{+})^{2}}(d_{+}(t)-d_{-}(t)).
\]
We would like to express the limit, as $\varepsilon_{-}+\varepsilon_{+}$ tends
to $0$, of the right-hand side term as a function of the limit quantities.
Remark that $u$ is intended to converge strongly but $\partial_{x}u$ only
weakly, thus $d_{+}(t)$ and $d_{-}(t)$ are not approximations of $\partial
_{x}u(t,x(t))$: however $\dfrac{\varepsilon_{-}}{\varepsilon_{-}%
+\varepsilon_{+}}d_{-}+\dfrac{\varepsilon_{+}}{\varepsilon_{-}+\varepsilon
_{+}}d_{-}$ is intended to converge toward $\partial_{x}u$. The limit of the
right-hand side should be expressed as a function of the limit unknowns
$\alpha_{+}$, $\alpha_{-}$, $p_{+}$, $p_{-}$, $\partial_{x}u$... We already
know that $\dfrac{\varepsilon_{-}\varepsilon_{+}}{(\varepsilon_{-}%
+\varepsilon_{+})^{2}}$ converges to $\alpha_{+}\alpha_{-}$. It remains to
treat the term $d_{+}-d_{-}$. As $\mu_{+}d_{+}-\mu_{-}d_{-}$ is intended to
converge to $p_{+}-p_{-}$, it is quite natural to try to write
\[
d_{+}-d_{-}=a(\mu_{+}d_{+}-\mu_{-}d_{-})+(1-a\mu_{+})d_{+}-(1-a\mu_{-})d_{-}%
\]
with $a\in\R$ such that there exists $b\in\R$ satisfying
\[
1-a\mu_{+}=b\alpha \quad\mbox{and}\quad1-a\mu_{-}=-b(1-\alpha),
\]
in which case one would have
\[
d_{+}-d_{-}\longrightarrow_{\varepsilon_{-}+\varepsilon_{+}\rightarrow
0}=a(p_{+}-p-)+b\partial_{x}u.
\]
The linear system in $a$ and $b$ has a unique solution, $a=\dfrac{1}%
{(1-\alpha)\mu_{+}+\alpha\mu_{-}}$ and $b=\dfrac{\mu_{-}-\mu_{+}}%
{(1-\alpha)\mu_{+}+\alpha\mu_{-}}$, which finally gives
\[
D_{t}\alpha_{+}=\frac{\alpha(1-\alpha)}{(1-\alpha)\mu_{+}+\alpha\mu_{-}}%
(p_{+}-p_{-}-(\mu_{+}-\mu_{-})\partial_{x}u),
\]
which is exactly the first equation in \eqref{TwoP}
\end{itemize}

\subsection{How to derive the effective viscosity and the effective pressure?} It is interesting to note
that the effective viscosity is computed in the same way as one computes the effective diffusion when considering homogenisation for 1d elliptic equations. The main remark is to use for the compressible the ellipticity of the effective flux $F= p(\rho) - \mu \partial_x u.$ Looking at the order $0$ terms in the momentum equation and averaging with respect to $y$ we get using that $u^0$ does not depend on $y$ that
$$\rho^0 (\partial_t u^0 + u^0 \partial_x u^0) - \partial_x (\overline{\mu^0}\partial_x u^0)
      - \partial_x (\overline{\mu^0 \partial_y u^1}) + \nabla \overline{p(\rho)} = 0.$$
 It remains to use the expression of $\overline{\mu^0 \partial_y u^1}$ found previously to get
 the homogeneized momentum equation related to $\alpha_\pm$ and $u^0$.

\section{Mathematical results\label{Mathematical_results}}

\subsection{Continuous approach -- Hoff solutions and its two-scale limit}

\noindent First of all, we can adapt the techniques introduced by D. Hoff in
order to obtain the well possedness of system \eqref{System}. More precisely,
let us assume the folllowing initial data conditions
\begin{equation}
\left\{
\begin{array}
[c]{l}%
c_{0}\in L^{\infty}({\mathbb{T}})\qquad\hbox{ such that }\qquad0\leq m=\inf
c_{0}\leq c_{0}(x)\leq M=\sup c_{0}\qquad\hbox{\it a.e. on }(0,1)\\
\rho_{0}\in L^{\infty}({\mathbb{T}})\qquad\hbox{ such that }\qquad0<\inf
\rho_{0}\leq\rho_{0}(x)\qquad\hbox{ a.e. on }(0,1)\\
G(\rho_{0},c_{0})=\rho_{0}%
{\displaystyle\int_{1}^{\rho_{0}}}
p(s,c_{0})/s^{2}<+\infty\\
u_{0}\in H^{1}({\mathbb{T}})
\end{array}
\right.  \label{INI}%
\end{equation}
and let us assume that $p$ and $\mu$ satisfies for all $(\rho,c)\in
\lbrack0,+\infty]\times\lbrack0,1]$ the following:
\begin{equation}
\left\{
\begin{array}
[c]{l}%
p(\rho,c)\geq0\quad\hbox{ and }\quad\mu(c)\geq\mu_{\mathrm{min}}>0\\
\exists C_{0}>0\hbox{ such that }p(\rho,c)\leq C_{0}(\rho+G(\rho
,c))\hbox{ where }G(\rho,c)=\rho\int_{1}^{\rho}p(\rho,c)/s^{2}\,ds\\
\rho\partial_{1}p(\rho,c)\in L_{\mathrm{loc}}^{\infty}([0,+\infty
]\times\lbrack0,1]).
\end{array}
\right.  \label{pc}%
\end{equation}

\bigskip

\noindent Then we have the following

\smallskip

\begin{theorem}
\label{Theorem1}\label{existence} Consider two functions $p\in\mathcal{C}%
^{1}\left(  [0,\infty)\times\left[  m,M\right]  \right)  $ and $\mu
\in\mathcal{C}^{1}\left(  \left[  0,M\right]  \right)  $ verifying Hypothesis
\eqref{pc}. Let $\left(  c_{0},\rho_{0},u_{0}\right)  $ satisfy \eqref{INI},
then there exists a unique weak solution $\left(  c,\rho,u\right)  $ of system
\eqref{System} with intial data $\left(  c_{0},\rho_{0},u_{0}\right)  $ with
\[%
\begin{array}
[c]{l}%
c,\rho\in C(\mathbb{[}0,\infty\mathbb{)};L^{q}\left(  \mathbb{T}^{1}\right)
)\hbox{ for all }q<+\infty,\\
u\in L^{\infty}\left(  \mathbb{[}0,\infty\mathbb{)};H^{1}(\mathbb{T}%
^{1})\right)  ,\text{ }\partial_{x}u\in L^{2}\left(  \mathbb{[}0,\infty
\mathbb{)}\times\mathbb{T}^{1}\right)  .
\end{array}
\]
Moreover, for any $T>0$, there exists a constant $C\left(  T\right)  $ which
depends only on the norms of the initial data and $T$ such that the following
uniform bounds hold true:%
\begin{gather}%
{\displaystyle\int_{0}^{1}}
\frac{\rho u^{2}}{2}+%
{\displaystyle\int_{0}^{1}}
G\left(  \rho,c\right)  +%
{\displaystyle\int_{0}^{t}}
{\displaystyle\int_{0}^{1}}
\mu(c)(\partial_{x}u)^{2}\leq%
{\displaystyle\int_{0}^{1}}
\frac{\rho_{0}u_{0}^{2}}{2}+%
{\displaystyle\int_{0}^{1}}
G\left(  \rho_{0},c_{0}\right)  ,\\
\inf\limits_{x\in\left[  0,1\right]  }c_{0}(x)\leq c(t,x)\leq\sup
\limits_{x\in\left[  0,1\right]  }c_{0}(x),\\
C\left(  T\right)  ^{-1}\leq\rho\left(  t,x\right)  \leq C\left(  T\right)
,\\
\frac{1}{2}%
{\displaystyle\int_{0}^{1}}
\mu(c)(\partial_{x}u)^{2}+%
{\displaystyle\int_{0}^{t}}
{\displaystyle\int_{0}^{1}}
\rho\dot{u}^{2}\leq C(T),\\
\frac{1}{2}%
{\displaystyle\int_{0}^{1}}
\sigma\left(  t\right)  \rho\dot{u}^{2}+\frac{1}{2}%
{\displaystyle\int_{0}^{t}}
{\displaystyle\int_{0}^{1}}
\sigma\left(  t\right)  \mu(c)\left(  \partial_{x}\dot{u}\right)  ^{2}\leq
C\left(  T\right)  ,\\
\left\Vert u\right\Vert _{L^{\infty}((0,1)\times\mathbb{T}^{1})}+\left\Vert
u\right\Vert _{H^{1}((0,1)\times\mathbb{T}^{1})}\leq C\left(  T\right)  ,\\
\sigma^{\frac{1}{2}}\left(  t\right)  \left\Vert (\partial_{t}u(t),\partial
_{x}u(t))\right\Vert _{L^{\infty}(\mathbb{T})}\leq C\left(  T\right)  .
\end{gather}
for all $t\in\left[  0,T\right]  $ where $\sigma\left(  t\right)
=\min\left\{  t,1\right\}  $.
\end{theorem}

\bigskip

\begin{remark}
It is interesting to note that the previous theorem includes general pressure
laws $p(\rho,c)$ and viscosity $\mu(c)$. It includes for instance the pressure
law $p(\rho,c)=\rho^{\gamma(c)}$ that depends on the density $\rho$ and the
fraction $c$ of each chemical/phase component. As described in \cite{MiVa},
the function $\gamma(c)$ depends on the constant heat capacity ratios of each
component of the multifluid, the pressure $p(\rho,c)$ effectively traces the
thermodynamic "signature" of mixing chemicals/phases in solution.
\end{remark}

\bigskip

Let us now consider a sequence of initial data $\left(  \rho_{0}^{\varepsilon
},c_{0}^{\varepsilon},u_{0}^{\varepsilon}\right)  _{\varepsilon>0}$ that
satisfies:
\begin{equation}
\left\{
\begin{array}
[c]{l}%
\rho_{0}^{\varepsilon}\in L^{\infty}\left(  \mathbb{T}^{1}\right)
\hbox{ with }0<\inf\limits_{n,x\in\mathbb{T}^{1}}\rho_{0}^{\varepsilon}\left(
x\right)  \leq\rho_{0}^{\varepsilon}\left(  x\right)  \leq\sup\limits_{n,x\in
\mathbb{T}^{1}}\rho_{0}^{\varepsilon}\left(  x\right)  \leq M<+\infty,\\
c_{0}^{\varepsilon}\,(1-c_{0}^{\varepsilon})=0\qquad\text{ a.e. on }%
\mathbb{T}^{1}\hbox{ with }\qquad c_{0}^{\varepsilon}\in\lbrack0,1],\\%
{\displaystyle\int_{0}^{1}}
\left(  c_{0}^{\varepsilon}p_{+}(\rho_{0}^{\varepsilon})+(1-c_{0}%
^{\varepsilon})p_{-}(\rho_{0}^{\varepsilon})\right)  \leq M,\\
u_{0}^{\varepsilon}\in H^{1}\left(  \mathbb{T}^{1}\right)  \text{ such that
}\left\Vert u_{0}^{\varepsilon}\right\Vert _{H^{1}}\leq M
\end{array}
\right.  \label{initial_data_H}%
\end{equation}
with $M>0$ independent of $\varepsilon$. We note that these assumptions are
satisfied in particular for initial configurations as depicted in
\eqref{eq_initialdata_scheme}. The bounds \eqref{initial_data_H} allow us to
conclude that there exists $\left(  \rho_{0},c_{0},u_{0}\right)  \in
L^{\infty}\left(  \mathbb{T}^{1}\right)  \times L^{\infty}\left(
\mathbb{T}^{1}\right)  \times H^{1}\left(  \mathbb{T}^{1}\right)  $ such that
\[
\rho_{0}^{\varepsilon}\rightharpoonup\rho_{0}\text{ in }L^{\infty}\left(
\mathbb{T}^{1}\right)  -w\star,\quad c_{0}^{\varepsilon}\rightharpoonup
\alpha_{0}\text{ in }L^{\infty}\left(  \mathbb{T}^{1}\right)  -w\star,\quad
u_{0}^{\varepsilon}\rightharpoonup u_{0}\text{ in }H^{1}\left(  \mathbb{T}%
^{1}\right)  .
\]
Furthermore, given $\varepsilon>0$ the initial data $(\rho_{0}^{\varepsilon
},c_{0}^{\varepsilon},u_{0}^{\varepsilon})_{\varepsilon>0}$ enters the scope
of Theorem \ref{Theorem1}. So, we can associate to this initial data a
solution $(\rho^{\varepsilon},c^{\varepsilon},u^{\varepsilon})_{\varepsilon
>0}$ to \eqref{initial_data_H}. Moreover, this sequence satisfies the
following uniform bounds on any interval $[0,T]$ independent of $\varepsilon$:%
\begin{gather}%
{\displaystyle\int_{0}^{1}}
\rho^{\varepsilon}(u^{\varepsilon})^{2}+%
{\displaystyle\int_{0}^{1}}
\left(  G(\rho^{\varepsilon},c^{\varepsilon})\right)  +\mu%
{\displaystyle\int_{0}^{t}}
{\displaystyle\int_{0}^{1}}
(\partial_{x}u^{\varepsilon})^{2}\leq C,\label{ineg(2)H}\\
C\left(  T\right)  ^{-1}\leq\rho^{\varepsilon}\left(  t,x\right)  \leq
C\left(  T\right)  ,\label{ineg(3)H}\\
\frac{\mu}{2}%
{\displaystyle\int_{0}^{1}}
(\partial_{x}u^{\varepsilon})^{2}+%
{\displaystyle\int_{0}^{t}}
{\displaystyle\int_{0}^{1}}
\rho^{\varepsilon}(\dot{u}^{\varepsilon})^{2}+\int_{0}^{t}\int_{0}%
^{1}\left\vert \partial_{x}(\mu\partial_{x}u^{\varepsilon}-p^{\varepsilon
})\right\vert ^{2}\leq C(T),\label{ineg(4)H}\\
\frac{1}{2}%
{\displaystyle\int_{0}^{1}}
\sigma\left(  t\right)  \rho^{\varepsilon}(\dot{u}^{\varepsilon})^{2}%
+\frac{\mu}{2}%
{\displaystyle\int_{0}^{t}}
{\displaystyle\int_{0}^{1}}
\sigma\left(  t\right)  \left(  \partial_{x}\dot{u}^{\varepsilon}\right)
^{2}\leq C\left(  T\right)  ,\label{ineg(5)H}\\
\left\Vert u^{\varepsilon}\right\Vert _{L^{\infty}((0,1)\times\mathbb{T}^{1}%
)}+\left\Vert u^{\varepsilon}\right\Vert _{H^{1}((0,1)\times\mathbb{T}^{1}%
)}\leq C\left(  T\right)  ,\label{ineg(6)H}\\
\sigma^{\frac{1}{2}}\left(  t\right)  \left\Vert (\partial_{t}u^{\varepsilon
}(t),\partial_{x}u^{\varepsilon}(t))\right\Vert _{L^{\infty}}\leq C\left(
T\right)  \label{ineg(7)}%
\end{gather}
with $\mu=\min(\mu_{+},\mu_{-}).$ Using the uniform bounds of $\left(
\text{\ref{ineg(2)H}}\right)  $-$\left(  \text{\ref{ineg(7)}}\right)  $ we
conclude that%
\begin{equation}
\left\{
\begin{array}
[c]{l}%
\rho{}^{\varepsilon}\rightharpoonup\rho{},\text{ }p{}(\rho{}^{\varepsilon
},c^{\varepsilon})\rightharpoonup\Pi{}\text{ in }L^{\infty}(\mathbb{R}%
_{+};L^{\infty}(\mathbb{T}^{1})),\\
u^{\varepsilon}\rightharpoonup u\text{ in }L^{\infty}(\mathbb{R}_{+}%
;H^{1}(\mathbb{T}^{1})),\\
Z^{n}:=\mu\partial_{x}u^{\varepsilon}-p(\rho^{\varepsilon},c^{\varepsilon
})\rightharpoonup Z^{\infty}:=\overline{\mu\partial_{x}u}-\Pi\text{ in }%
L^{2}(\mathbb{R}_{+};H^{1}(\mathbb{T}^{1})).
\end{array}
\right.  \label{weak_limits}%
\end{equation}
As explained previously, the density $\rho^{\varepsilon}$ and the parameter
$c^{\varepsilon}$ are expected to oscillate widely in space. For this reason,
it is hopeless to obtain stronger convergence on these sequences (and it would
not be a good news for our plan either) than in a weak $L^{p}$-setting. It
turns out that we may characterize weak*-limits of the form if we suppose that
they are given at time $0$. In the following lines we make rigourous the
previous assertion. On the other hand, we need to recover some properties of
the sequence $p(\rho^{\varepsilon},c^{\varepsilon})$ to compute a limit system
satisfied by $(\rho,u,\Pi).$ To this end, we associate to the sequence
$(\rho^{\varepsilon},c^{\varepsilon})_{\varepsilon}$ a sequence of measures on
the space $\mathbb{T}_{x}^{1}\times\mathbb{R}_{\xi}\times\mathbb{R}_{\eta}$
(here $\mathbb{R}_{\xi}$ must be understood as the range of the $\rho
^{\varepsilon}$ while $\mathbb{R}_{\eta}$ is the range of the $c^{\varepsilon
}$). Namely, given $n\geq0$ and $t\geq0$, we consider the measure on
$\mathbb{T}_{x}^{1}\times\mathbb{R}_{\xi}\times\mathbb{R}_{\eta}$ as defined
by%
\begin{equation}
\left\langle \Theta{}^{\varepsilon}\left(  t\right)  ,b\right\rangle
:\overset{def.}{=}\int_{\mathbb{T}^{1}}b\left(  x,\rho{}^{\varepsilon}\left(
t,x\right)  ,c^{\varepsilon}(t,x)\right)  dx,\qquad\forall\,b\in C_{c}\left(
\mathbb{T}_{x}^{1}\mathbb{\times R}_{\xi}\times\mathbb{R}_{\eta}\right)
\label{suite_mesures}%
\end{equation}

\bigskip

\noindent We have the following proposition.

\begin{proposition}
For fixed $\varepsilon>0$ fixed there holds
\begin{equation}
\Theta^{\varepsilon}\in C_{w}([0,\infty);\mathcal{M_{+}(\mathbb{T}}_{x}%
^{1}\mathcal{\mathbb{\times R}_{\xi}\times\mathbb{R}_{\eta})})
\label{eq_continuite}%
\end{equation}
with
\begin{equation}
\operatorname{Supp}(\Theta^{\varepsilon}(t))\subset\mathbb{T}_{x}^{1}%
\times\lbrack C\left(  t\right)  ^{-1},C\left(  t\right)  ]\times
\lbrack0,1]\quad\langle\Theta^{\varepsilon},1\rangle=1.\quad\forall\,t\geq0,
\label{eq_supportThetan}%
\end{equation}
where $C\left(  t\right)  $ is given by \eqref{ineg(3)H}.
\end{proposition}

\bigskip

\noindent Once these measures are constructed, the rigorous justification of
system \eqref{TwoP} consists in the derivation of the kinetic limit and the characteirzation
of the measures family, it  reduces to

\begin{theorem}
\label{Theorem3} Up to the extraction of a subsequence, we have $\Theta
^{\varepsilon}\rightharpoonup\Theta${ in $C_{w}([0,\infty);\mathcal{M_{+}%
(}\mathbb{T}_{x}^{1}\times\mathbb{R}_{\xi}\times\mathbb{R}_{\eta}))$} where
$\Theta$ satisfies
\begin{equation}
\partial_{t}\Theta+\partial_{x}\left(  u\,\Theta{}\right)  -\partial_{\xi
}\left(  \left(  \dfrac{\xi Z^{\infty}}{\mu{}\left(  \eta\right)  }+\dfrac{\xi
p{}(\xi,\eta)}{\mu{}\left(  \eta\right)  }\right)  \Theta\right)  -\left(
\dfrac{Z^{\infty}}{\mu\left(  \eta\right)  }+\dfrac{p{}\left(  \xi
,\eta\right)  }{\mu{}\left(  \eta\right)  }\right)  \Theta
=0\label{equation_of_thetalim}%
\end{equation}
with $(u,\Pi,Z^{\infty})$ as defined in \eqref{weak_limits}. Moreover, if
there exists $(\alpha_{0},\rho_{+,0},\rho_{-,0})\in L^{\infty}(\mathbb{T}%
^{1})$ with $\alpha_{0}\in\left[  0,1\right]  $ a.e. and such that
\begin{multline}
\langle\Theta(0),b\rangle=\int_{\mathbb{T}^{1}}(\alpha_{0}(x)b\left(
x,\rho_{+,0}(x),0\right)  +(1-\alpha_{0}(x))b\left(  x,\rho_{-,0}(x),1\right)
)dx,\\
\quad\forall\,b\in C(\mathbb{T}_{x}^{1}\times\mathbb{R}_{\xi}\times
\mathbb{R}_{\eta})\label{struct_initiala}%
\end{multline}
then there exists $(\alpha,\rho_{+},\rho_{-})\in\lbrack L_{loc}^{\infty
}([0,\infty);L^{\infty}(\mathbb{T}))\cap C([0,\infty);L^{1}(\mathbb{T}))]^{4}$
such that, for any $t\geq0,$ $\alpha\left(  t\right)  \in\left[  0,1\right]  $
a.e. and
\begin{multline}
\langle\Theta(t),b\rangle=\int_{\mathbb{T}^{1}}(\alpha(t,x)b\left(
x,\rho_{+,0}(t,x),0\right)  +(1-\alpha(t,x))b\left(  x,\rho_{-,0}%
(t,x),1\right)  )dx,\\
\quad\forall\,b\in C(\mathbb{T}_{x}^{1}\times\mathbb{R}_{\xi}\times
\mathbb{R}_{\eta}),\label{structure}%
\end{multline}
Furthermore, $\left(  \alpha,\rho_{+},\rho_{-}\right)  $ together with $u$
verifies the biphase Baer-Nunziato type system \eqref{TwoP}.
\end{theorem}

\noindent {\bf Remark.} It is interesting to understand that Hoff's solution
allows the density to be only bounded and the velocity to be sufficiently 
regular to characterize its evolution in time.

\subsection{The semi-discrete approach -- ODEs and its continuous limit}

This part corresponds in some sense to a mathematical justification of the
formal description indicated previous. We consider the following system of
ODEs:%
\begin{equation}
\left\{
\begin{array}
[c]{l}%
\dot{x}_{j+\frac{1}{2}}=u_{j+\frac{1}{2}},\\
\dot{c}_{j}=0,\\
\dfrac{d}{dt}\left(  \rho_{j}\Delta x_{j}\right)  =0,\\
\rho_{j+\frac{1}{2}}\Delta x_{j+\frac{1}{2}}\dot{u}_{j+\frac{1}{2}}%
+p_{j+1}-p_{j}=\left\{  \mu\left(  c_{j+1}\right)  \dfrac{u_{j+\frac{3}{2}%
}-u_{j+\frac{1}{2}}}{\Delta x_{j+1}}-\mu\left(  c_{j}\right)  \dfrac
{u_{j+\frac{1}{2}}-u_{j-\frac{1}{2}}}{\Delta x_{j}}\right\}  ,\\
\Delta x_{j}=x_{j+\frac{1}{2}}-x_{j-\frac{1}{2}},\text{ \ }\\
\Delta x_{j+\frac{1}{2}}=\dfrac{\Delta x_{j}+\Delta x_{j+1}}{2},\\
\rho_{j+\frac{1}{2}}=\dfrac{\rho_{j}\Delta x_{j}+\rho_{j+1}\Delta x_{j+1}%
}{\Delta x_{j}+\Delta x_{j+1}}.
\end{array}
\right.  \label{semi_discrete_semi_lagrangian}%
\end{equation}
for all $j\in\overline{0,J-1}$ and
\begin{equation}
\left\{
\begin{array}
[c]{c}%
c_{0}=c_{J},\\
\rho_{0}=\rho_{J},\\
u_{\frac{1}{2}}=u_{J+\frac{1}{2}}.
\end{array}
\right.  \label{periodic_boundary_condition}%
\end{equation}
From the previous system of equations we also deduce that%
\begin{equation}
\left\{
\begin{array}
[c]{c}%
\dot{\Delta}x_{j}=u_{j+\frac{1}{2}}-u_{j-\frac{1}{2}},\\
\dfrac{d}{dt}\left(  \rho_{j+\frac{1}{2}}\Delta x_{j+\frac{1}{2}}\right)  =0.
\end{array}
\right.  \label{auxilary}%
\end{equation}

\noindent System $\left(  \text{\ref{semi_discrete_semi_lagrangian}}\right)  $
is to be completed with initial data $\left(  x_{j+\frac{1}{2}}^{0}\right)
_{j=\overline{-1,J-1}}$, $\left(  c_{j}^{0}\right)  _{j=\overline{0,J-1}}$,
$\left(  \rho_{j}^{0}\right)  _{j=\overline{0,J-1}}$, $\left(  u_{j+\frac
{1}{2}}^{0}\right)  _{j=\overline{0,J-1}}$ such that:%
\begin{equation}
\left\{
\begin{array}
[c]{l}%
x_{-\frac{1}{2}}^{0}<\text{ }x_{\frac{1}{2}}^{0}<x_{\frac{3}{2}}^{0}%
<\cdots<x_{J-\frac{1}{2}}^{0},\\
0<\underline{c^{0}}=\min\limits_{j\in\overline{0,J-1}}c_{j}^{0}\leq
\overline{c^{0}}=\max\limits_{j\in\overline{0,J-1}}c_{j}^{0}<\infty,\\
0<\underline{\rho^{0}}=\min\limits_{j\in\overline{0,J-1}}\rho_{j}^{0}%
\leq\overline{\rho^{0}}=\max\limits_{j\in\overline{0,J-1}}\rho_{j}^{0}%
<\infty,\\
\left\Vert \left(  u_{j+\frac{1}{2}}\right)  _{j\in\overline{0,J-1}%
}\right\Vert _{\hat{H}_{J}^{1}}^{2}=%
{\displaystyle\sum\limits_{j=0}^{J-1}}
\left\vert u_{j+\frac{1}{2}}^{0}\right\vert ^{2}\Delta x_{j}^{0}+%
{\displaystyle\sum\limits_{j=0}^{J-1}}
\left\vert \dfrac{u_{j+\frac{1}{2}}^{0}-u_{j-\frac{1}{2}}^{0}}{\Delta x_{j}%
}\right\vert ^{2}\Delta x_{j}^{0}<\infty.
\end{array}
\right.  \tag{$\textrm{H}$}\label{discrete_Initial_data}%
\end{equation}

\bigskip

\noindent System $\left(  \text{\ref{semi_discrete_semi_lagrangian}}\right)  $
is a system of ODEs which, owing to the fact that%
\[
c_{j}\left(  t\right)  =c_{j}\left(  0\right)  ,\rho_{j}\left(  t\right)
=\dfrac{\rho_{j}\left(  0\right)  \Delta x_{j}\left(  0\right)  }%
{x_{j+\frac{1}{2}}\left(  t\right)  -x_{j-\frac{1}{2}}\left(  t\right)  }%
,\rho_{j+\frac{1}{2}}\left(  t\right)  \Delta x_{j+\frac{1}{2}}\left(
t\right)  =\rho_{j+\frac{1}{2}}\left(  0\right)  \Delta x_{j+\frac{1}{2}%
}\left(  0\right)  ,
\]
can be put in the form%
\begin{equation}
\left\{
\begin{array}
[c]{l}%
\dfrac{dX}{dt}=U,\\
\dfrac{dU}{dt}=F\left(  X,U\right)  ,\\
\left(  X,U\right)  =\left(  X_{0},U_{0}\right)  ,
\end{array}
\right.  \label{system_ode}%
\end{equation}
with $X=\left(  x_{j+\frac{1}{2}}\right)  _{j=0,J-1}\in\mathbb{R}^{J}$ and
$U=\left(  u_{j+\frac{1}{2}}\right)  _{j=0,J-1}\in\mathbb{R}^{J}$ and
$F:D\rightarrow\mathbb{R}^{J}$ where
\[
D=\left\{  \left(  X,U\right)  \in\mathbb{R}^{J}\times\mathbb{R}^{J}:\text{
}x_{\frac{1}{2}}<x_{\frac{3}{2}}<\cdots<x_{J-\frac{1}{2}}\right\}
\]
which is an open set of $\mathbb{R}^{J}\times\mathbb{R}^{J}$. Owing to the
fact that $F$ is $C^{\infty}$ on $D$ we obtain via the Cauchy-Lipschitz-Peano
theorem that for any initial data $\left(  X_{0},U_{0}\right)  \in D$, there
exists a unique maximal solution for $\left(  \text{\ref{system_ode}}\right)
$
\[
\left(  X,U\right)  :[0,T_{\max})\rightarrow D
\]
with $T_{\max}=\infty$ or if $T_{\max}<\infty$ then
\begin{equation}
\lim_{t\rightarrow T_{\max}}\left(  X\left(  t\right)  ,U\left(  t\right)
\right)  \in\partial D \label{explosion}%
\end{equation}

As expected, system $\left(  \text{\ref{semi_discrete_semi_lagrangian}%
}\right)  $ shares a lot of properties with its continious version, in
particular we can find a number of apriori estimates related to the mass and
energy conservation which enable us to prove that $\left(
\text{\ref{explosion}}\right)  $ does not occur. We are thus able to show that
the local solutions can be extended to global ones. The following energy
functionals play a key role: the total mass Consider the mass%
\begin{equation}
M\left(  t\right)  =\text{\ }%
{\displaystyle\sum\limits_{j=0}^{J-1}}
\rho_{j}\left(  t\right)  \Delta x_{j}\left(  t\right)  \label{mass}%
\end{equation}
and basic energy functionals
\begin{multline}
E\left(  t\right)  =\dfrac{1}{2}%
{\displaystyle\sum\limits_{j=0}^{J-1}}
\rho_{j+\frac{1}{2}}\left(  t\right)  \left\vert u_{j+\frac{1}{2}}\right\vert
^{2}\Delta x_{j+\frac{1}{2}}\left(  t\right)  +%
{\displaystyle\sum\limits_{j=0}^{J-1}}
H(\rho_{j}\left(  t\right)  ,c_{j}\left(  t\right)  )\Delta x_{j}\left(
t\right) \\
+%
{\displaystyle\int\limits_{0}^{t}}
{\displaystyle\sum\limits_{j=0}^{J-1}}
\mu\left(  c_{j}\left(  \tau\right)  \right)  \left\vert \dfrac{u_{j+\frac
{1}{2}}\left(  \tau\right)  -u_{j-\frac{1}{2}}\left(  \tau\right)  }{\Delta
x_{j}\left(  \tau\right)  }\right\vert ^{2}\Delta x_{j}\left(  \tau\right)  .
\label{energy}%
\end{multline}

\bigskip

\noindent We formalize in the following theorem our first result.

\begin{theorem}
\label{Numeric}The system \eqref{semi_discrete_semi_lagrangian} along with
initial data verifying hypothesis $\left(  \text{\ref{discrete_Initial_data}%
}\right)  $ admits an unique solution and, moreover, the following estimates
hold true uniformly in $J$:%
\[
\left\{
\begin{array}
[c]{l}%
\dfrac{1}{C_{ini}^{1}\left(  t\right)  }\leq\rho_{j}\left(  t\right)  \leq
C_{ini}^{1}\left(  t\right)  ,\text{ for all }j\in\overline{0,J-1},\\
\dfrac{1}{C_{ini}^{2}\left(  t\right)  }\Delta x_{j}^{0}\leq\Delta
x_{j}\left(  t\right)  \leq\Delta x_{j}^{0}C_{ini}^{2}\left(  t\right)
,\text{ for all }j\in\overline{0,J-1},\\
\underline{c^{0}}\leq c_{j}\left(  t\right)  \leq\overline{\text{ }c^{0}}\\
M\left(  t\right)  =M_{0},\text{ }E\left(  t\right)  =E_{0}%
\end{array}
\right.
\]
where $C_{ini}^{1}\left(  \cdot\right)  $, $C_{ini}^{2}\left(  \cdot\right)  $
and $C_{ini}^{3}\left(  \cdot\right)  $ are strictly positive increasing
continuous functions that depend only on the initial data.
\end{theorem}

\medskip

For any $J\in\mathbb{N}^{\ast}$, having constructed the functions $\left(
c_{j},\rho_{j},u_{j+\frac{1}{2}}\right)  _{j\in\overline{0,J-1}}$ as above we
consider
\[
\hat{c}_{J},\hat{\rho}_{J},\hat{u}_{J}:[0,\infty)\times\mathbb{R\rightarrow R}%
\]
defined by%
\begin{align}
\hat{c}_{J}\left(  t,x\right)   &  =c_{j}\left(  t\right)  \text{ if }%
x\in\lbrack x_{j-\frac{1}{2}}\left(  t\right)  ,x_{j+\frac{1}{2}}\left(
t\right)  )\text{,}\label{def_c_cont}\\
\hat{\rho}_{J}\left(  t,x\right)   &  =\rho_{j}\left(  t\right)  \text{ if
}x\in\lbrack x_{j-\frac{1}{2}}\left(  t\right)  ,x_{j+\frac{1}{2}}\left(
t\right)  )\text{,}\label{def_rho_cont}\\
\hat{u}_{J}\left(  t,x\right)   &  =\frac{x-x_{j-\frac{1}{2}}}{\Delta x_{j}%
}u_{j+\frac{1}{2}}\left(  t\right)  +\frac{x_{j+\frac{1}{2}}-x}{\Delta x_{j}%
}u_{j-\frac{1}{2}}\left(  t\right)  \text{ if }x\in\lbrack x_{j-\frac{1}{2}%
}\left(  t\right)  ,x_{j+\frac{1}{2}}\left(  t\right)  )\text{,}%
\label{def_u_cont}\\
\hat{Z}_{J}\left(  t,x\right)   &  =\frac{x-x_{j}}{\Delta x_{j+\frac{1}{2}}%
}\sigma_{j+1}+\frac{x_{j+1}-x}{\Delta x_{j+\frac{1}{2}}}\sigma_{j}\text{ if
}x\in\lbrack x_{j},x_{j+1}).\label{effective_flux}%
\end{align}
where%
\[
Z_{j}=\mu\left(  c_{j}\right)  \frac{u_{j+\frac{1}{2}}-u_{j-\frac{1}{2}}%
}{\Delta x_{j}}-p\left(  c_{j},\rho_{j}\right)  .
\]

\medskip

\noindent First, we observe the following remarkable property which will be
strongly use to get the limit

\begin{proposition}
\label{transport_discret}The functions $\left(  \hat{c}_{J},\hat{\rho}%
_{J},\hat{u}_{J}\right)  $ verify the following transport equations
\begin{equation}
\left\{
\begin{array}
[c]{l}%
\partial_{t}\hat{c}_{J}+\hat{u}_{J}\partial_{x}\hat{c}_{J}=0,\\
\partial_{t}\hat{\rho}_{J}+\partial_{x}\left(  \hat{\rho}_{J}\hat{u}%
_{J}\right)  =0,
\end{array}
\right.  \label{transport_discrete}%
\end{equation}
with initial data
\[
\left\{
\begin{array}
[c]{l}%
\hat{\rho}_{J|t=0}=\hat{\rho}_{J}^{0},\\
\hat{c}_{J|t=0}=\hat{c}_{J}^{0},
\end{array}
\right.
\]
in the sense of distributions.
\end{proposition}

\bigskip

\noindent Of course, the estimates announced in Theorem \ref{Numeric} can be
used in order to estimate various norms of the functions $\left(  \hat{c}%
_{J},\hat{\rho}_{J},\hat{u}_{J},\hat{Z}_{J}\right) $ as in the Hoff solution proof.

\medskip

\noindent More precisely, one has that

\begin{theorem}
\label{Thm_cont}Consider discrete initial data verifying the hypothesis
\eqref{discrete_Initial_data} along with the globally defined solution of the
system of ODEs
\eqref{semi_discrete_semi_lagrangian}--\eqref{periodic_boundary_condition}.
Furthermore, consider the functions $\left(  \hat{c}_{J},\hat{\rho}_{J}%
,\hat{u}_{J},\hat{Z}_{J}\right)  $ given by
\eqref{def_c_cont}--\eqref{effective_flux}. Then
\begin{equation}
\left\{
\begin{array}
[c]{l}%
\min\limits_{j\in\overline{0,J-1}}c_{j}^{0}\leq\hat{c}_{J}\left(  t,x\right)
\leq\max\limits_{j\in\overline{0,J-1}}c_{j}^{0},\\
\dfrac{1}{C_{ini}^{1}\left(  t\right)  }\leq\hat{\rho}_{J}\left(  t,x\right)
\leq C_{ini}^{1}\left(  t\right)  ,\\%
{\displaystyle\int_{0}^{1}}
\hat{\rho}_{J}\left\vert \hat{u}_{J}\right\vert ^{2}\left(  t,x\right)  dx+%
{\displaystyle\int_{0}^{1}}
H\left(  \hat{\rho}_{J}\left(  t,x\right)  ,\hat{c}_{J}\left(  t,x\right)
\right)  dx+%
{\displaystyle\int_{0}^{t}}
{\displaystyle\int_{0}^{1}}
\mu\left(  \hat{c}_{J}\left(  \tau,x\right)  \right)  \left\vert \partial
_{x}\hat{u}_{J}\left(  \tau,x\right)  \right\vert ^{2}dxd\tau\leq2E_{0},\\
\left\Vert \partial_{x}\hat{u}_{J}\right\Vert _{L_{t,x}^{2}}+\left\Vert
\partial_{t}\hat{u}_{J}\right\Vert _{L_{t,x}^{2}}+\min\left\{  1,t\right\}
\left\Vert \partial_{t}\hat{u}_{J}\left(  t\right)  \right\Vert _{L_{x}^{2}%
}\leq C\left(  t\right)  ,\\%
{\displaystyle\int_{0}^{t}}
{\displaystyle\int_{0}^{1}}
\left\vert \partial_{x}\hat{\sigma}_{J}\left(  \tau,x\right)  \right\vert
^{2}dxd\tau+\min\left\{  1,t\right\}
{\displaystyle\int_{0}^{1}}
\left\vert \partial_{x}\hat{\sigma}_{J}\left(  t,x\right)  \right\vert ^{2}dx+%
{\displaystyle\int_{0}^{t}}
(\sup_{x\in\left[  0,1\right]  }\left\vert \partial_{x}\hat{\sigma}_{J}\left(
\tau,x\right)  \right\vert )^{\frac{4}{3}-}d\tau\leq C\left(  t\right) \\%
{\displaystyle\int_{0}^{t}}
{\displaystyle\int_{0}^{1}}
\sigma\left(  \tau\right)  \left\vert \partial_{t}\hat{\sigma}_{J}\left(
\tau,x\right)  \right\vert ^{2}dxd\tau\leq C\left(  t\right)  .
\end{array}
\right.  \label{Estimari_cont}%
\end{equation}

\end{theorem}

and basic energy functionals

\bigskip

For all $J\in\mathbb{N}^{\ast}$ and all $t\geq0$, we consider the measure on
$\mathbb{T}_{x}^{1}\times\mathbb{R}_{\xi}\times\mathbb{R}_{\eta}$ defined by%
\begin{equation}
\left\langle \Theta{}^{J}\left(  t\right)  ,b\right\rangle :\overset{def.}%
{=}\int_{\mathbb{T}^{1}}b\left(  x,\hat{\rho}{}^{J}\left(  t,x\right)
,\hat{c}^{J}(t,x)\right)  dx,\qquad\forall\,b\in C_{c}\left(  \mathbb{T}%
_{x}^{1}\mathbb{\times R}_{\xi}\times\mathbb{R}_{\eta}\right)
\label{measure_disc}%
\end{equation}
We may establish immediately similar results as in the previous section namely:

\begin{proposition}
For fixed $J\in\mathbb{N}$ we have that
\begin{equation}
\Theta^{J}\in C_{w}([0,\infty);\mathcal{M_{+}(\mathbb{T}}_{x}^{1}%
\mathcal{\mathbb{\times R}_{\xi}\times\mathbb{R}_{\eta})})
\label{eq_continuite_disc}%
\end{equation}
with
\begin{equation}
\operatorname{Supp}(\Theta^{J}(t))\subset\mathbb{T}_{x}^{1}\times\lbrack
C_{ini}^{1}\left(  t\right)  ^{-1},C_{ini}^{1}\left(  t\right)  ]\times
\lbrack0,1]\quad\langle\Theta^{J},1\rangle=1.\quad\forall\,t\geq0,
\label{eq_supportThetan_disc}%
\end{equation}
where $C_{ini}^{1}\left(  t\right)  $ is given by
\eqref{Estimari_cont}.\bigskip
\end{proposition}

We can also pass to the limit in the kinetic system to get the following result

\begin{theorem}
\label{Theorem3_discret} Up to the extraction of a subsequence, we have
$\Theta^{J}\rightharpoonup\Theta${ in $C_{w}([0,\infty);\mathcal{M_{+}%
(}\mathbb{T}_{x}^{1}\times\mathbb{R}_{\xi}\times\mathbb{R}_{\eta}))$} where
$\Theta$ satisfies
\begin{equation}
\partial_{t}\Theta+\partial_{x}\left(  u\,\Theta{}\right)  -\partial_{\xi
}\left(  \left(  \dfrac{\xi Z^{\infty}}{\mu{}\left(  \eta\right)  }+\dfrac{\xi
p{}(\xi,\eta)}{\mu{}\left(  \eta\right)  }\right)  \Theta\right)  -\left(
\dfrac{Z^{\infty}}{\mu\left(  \eta\right)  }+\dfrac{p{}\left(  \xi
,\eta\right)  }{\mu{}\left(  \eta\right)  }\right)  \Theta
=0\label{equation_of_thetalim_disc}%
\end{equation}
with $(u,\Pi,Z)$ as defined in \eqref{weak_limits}. Moreover, if there exists
$(\alpha_{0},\rho_{+,0},\rho_{-,0})\in L^{\infty}(\mathbb{T}^{1})$ such that
$\alpha_{0}\in\left[  0,1\right]  $ a.e. and
\begin{multline}
\langle\Theta(0),b\rangle=\int_{\mathbb{T}^{1}}(\alpha_{0}(x)b\left(
x,\rho_{+,0}(x),0\right)  +(1-\alpha_{0}(x))b\left(  x,\rho_{-,0}(x),1\right)
)dx\\
\quad\forall\,b\in C(\mathbb{T}_{x}^{1}\times\mathbb{R}_{\xi}\times
\mathbb{R}_{\eta}),\label{struct_initiala_disc}%
\end{multline}
then there exists $(\alpha,\rho_{+},\rho_{-})\in\lbrack L_{loc}^{\infty
}([0,\infty);L^{\infty}(\mathbb{T}))\cap C([0,\infty);L^{1}(\mathbb{T}))]^{4}$
such that, for any $t\geq0$ we have $\alpha\left(  t,\cdot\right)  \in\left[
0,1\right]  $ a.e. and
\begin{multline}
\langle\Theta(t),b\rangle=\int_{\mathbb{T}^{1}}(\alpha(t,x)b\left(  x,\rho
_{+}(t,x),0\right)  +(1-\alpha(t,x))b\left(  x,\rho_{-}(t,x),1\right)  )dx\\
\quad\forall\,b\in C(\mathbb{T}_{x}^{1}\times\mathbb{R}_{\xi}\times
\mathbb{R}_{\eta}).\label{structure_disc}%
\end{multline}
Furthermore, $\left(  \alpha,\rho_{+},\rho_{-}\right)  $ together with $u$
verifies the biphase Baer-Nunziato system \eqref{TwoP}.
\end{theorem}

\section{Numerical illustrations}\label{Numerical_illustrations}

In this last section we will illustrate the results stated in the paper.
We design two numerical schemes: one to approximate the mesoscopic system, that is to say System \eqref{meso} with a Cauchy datum \eqref{mesoCauchy}, and one to approximate the macroscopic system \eqref{TwoP} with any Cauchy datum.

\subsection{Mesoscopic discretization}\label{mesoscheme}

The numerical scheme we design here consists in a ``brute force"
discretization of
System~\eqref{meso} where $c_0(1-c_0) = 0$.

As the fluids have to remain pure (not mixed) in every cell, because,
for modelling reasons, we want to use
only the pure pressure laws (the mixture pressure law being unknown
at this stage), the length of each
pure zone has to be larger than a cell (and, more precisely, has to
be large as an integer number of cells).
Here, in the numerical tests, we choose to consider a numerical
initial condition such that the fluid changes
from one cell to the other (but of course this is not a restriction).
The problem to achieve the aim
here comes from the so-called numerical diffusion:
the discretization of the $\partial_t c + u \partial_x c = 0$
with a stable scheme usually brings a certain amount of diffusion, the
effect of which being not to preserve the important feature
$c(t,\cdot)(1-c(t,\cdot)) = 0 \quad a.e.$.
In order to pass over this phenomenon, we consider a {\em Lagrangian},
or {\em pseudo-Lagrangian}\footnote{It can be called pseudo-Lagrangian because,
although the solution is actually expressed in the classical Euler variable,
the scheme strongly uses the Langrange formulation of the system. }
scheme in which the cells follow the fluid in its transport,
namely in which the edges of every cell moves at the fluid
velocity. In this Lagrangian frame, the equation for the mass fraction
is $D_t c = 0$ (recall that $D_t = \partial_t + u \partial_x$).

The spirit of the proposed scheme is the one of staggered schemes:
it can be seen as a modification of the schemes in \cite{Karper} and
\cite{HKL}, this modification being that the present scheme is more
explicit (precisely, the nonlinearity are time-discretized in a backward
Euler way) and that it is a pseudo-Lagrange scheme. Staggered schemes are schemes in which different unknowns are associated to different points or cells in the mesh (for example, the density and the velocity, here).
At last, this scheme is a time discretization of the semi-discrete scheme \eqref{semi_discrete_semi_lagrangian} that was proposed to determine the limit macroscopic system.


The discretization is the following. Let $J \in \N \setminus \{0\}$
be the number of cells in $[0, 1)$. Let
$(x_{j-1/2}^0)_{j = 1}^{J}$ be the collection of cell interface positions
at time $0$. One assumes
$0 \leq x_{j-1/2}^0 < x_{j+1/2}^0 < 1$ for any $j = 1,\dots,J-1$.

In order to take into account the fact that the problem under
consideration is posed on $\T$ in a simple manner, i.e. without taking care
of the cells and quantities on the boundary, we extend all the data over $\R$ and $\Z$ by periodicity.

The cells itselves are denoted by $\omega_j^0 = [x_{j-1/2}^0, x_{j+1/2}^0)$
for $j \in \Z$. We denote by $\Delta x_j^n = x_{j+1/2}^n - x_{j-1/2}^n$ their length. The minimum length of
these cells is intended
to be small (and to tend to $0$ as $J$ tends to $\infty$ to reach convergence).
We also will need the distance between two centers of consecutive cells:
$\Delta x_{j+1/2}^n = (\Delta x_j^n + \Delta x_{j+1}^n)/2$.

Each time step of the scheme, given a discrete datum
$\left(x_{j-1/2}^n,\rho_j^n,c_j^n,u_{j-1/2}^n\right)_{j \in \Z}$,
consists in defining appropriately $\Delta t^n > 0$
and constructing $\left(x_{j-1/2}^{n+1},\rho_j^{n+1},c_j^{n+1},u_{j-1/2}^{n+1}\right)_{j \in \Z}$ by the formula
\begin{equation}\label{hypsch}
\left\{
\begin{array}{l}
\displaystyle{}\rho_{j+1/2}^n = \frac{\Delta x_j^n \rho_{j}^n + \Delta x_{j+1}^n \rho_{j+1}^n}{\Delta x_{j}^n + \Delta x_{j+1}^n}, \quad j \in \Z, \\
\displaystyle{} c_j^{n+1} = c_j^n, \quad j \in \Z, \\
\displaystyle{} \rho_{j+1/2}^{n} \Delta x_{j+1/2}^{n} u_{j+1/2}^{n+1} =
\rho_{j+1/2}^{n} \Delta x_{j+1/2}^{n} u_{j+1/2}^{n} - \Delta t^n
\left( p(c_{j+1}^n,\rho_{j+1}^n) - p(c_{j}^n,\rho_{j}^n) \right) \\
\displaystyle{}\phantom{\displaystyle{} \rho_{j+1/2}^{n+1} \Delta x_{j+1/2}^{n+1} u_{j+1/2}^{n+1} = } + \Delta t^n
\left( \mu(c_{j+1}^n,\rho_{j+1}^n) \frac{u_{j+3/2}^{n+1} - u_{j+1/2}^{n+1}}{\Delta x_{j+1}^n}
 - \mu(c_{j}^n,\rho_j^n)\frac{u_{j+1/2}^{n+1} - u_{j-1/2}^{n+1}}{\Delta x_{j}^n} \right), \quad j \in \Z, \\
x_{j+1/2}^{n+1} = x_{j+1/2}^{n} + \Delta t^n u_{j+1/2}^{n+1}, \quad j \in \Z, \\
\displaystyle{} \Delta x_j^{n+1} = x_{j+1/2}^{n+1} - x_{j-1/2}^{n+1},
\quad j \in \Z, \\
\displaystyle{} \Delta x_{j+1/2}^{n+1} = \frac{\Delta x_{j}^{n+1}
+ \Delta x_{j+1}^{n+1}}{2}, \quad j \in \Z, \\
\displaystyle{} \rho_{j}^{n+1} = \rho_j^n \frac{\Delta x_j^n}{\Delta x_j^{n+1}},
\quad j \in \Z. \\
\end{array}
\right.
\end{equation}
In the system above,
\begin{itemize}
\item The first equation defines a density associated to the nodes $x_{j+1/2}^n$, density that is used in the third equation,
\item The second equation is a (non-diffusive) discretization of $D_t c = 0$,
\item The third equation is the discretization of $\partial_t \rho u + \partial_x (\rho u^2 + p) = \partial_x(\mu\partial_x u)$: indeed notice that thanks to the last equation of the system, this third equation rewrites
\begin{multline*}
\displaystyle{} \rho_{j+1/2}^{n+1} \Delta x_{j+1/2}^{n+1} u_{j+1/2}^{n+1} =
\rho_{j+1/2}^{n} \Delta x_{j+1/2}^{n} u_{j+1/2}^{n} - \Delta t^n
\left( p(c_{j+1}^n,\rho_{j+1}^n) - p(c_{j}^n,\rho_{j}^n) \right) \\
\displaystyle{} + \Delta t^n
\left( \mu(c_{j+1}^n,\rho_{j+1}^n) \frac{u_{j+3/2}^{n+1} - u_{j+1/2}^{n+1}}{\Delta x_{j+1}^n}
 - \mu(c_{j}^n,\rho_j^n)\frac{u_{j+1/2}^{n+1} - u_{j-1/2}^{n+1}}{\Delta x_{j}^n} \right), \quad j \in \Z,
\end{multline*}
which is consistent with the partial differential equation,
\item The fourth equation is the translation of the mesh,
\item Fifth and sixth equations redefine quantities that are used in the scheme,
\item The last equation expresses the conservation of mass in a material volume $\partial_t \rho + \partial_x \rho u = 0$).
\end{itemize}
It is possible to prove that if the time step $\Delta t^n$ is sufficently small,
$x_{j-1/2}^n < x_{j+1/2}^n$ for all $j$ implies $x_{j-1/2}^{n+1} < x_{j+1/2}^{n+1}$ for all $j$.

\subsection{Macroscopic discretization}\label{macroscheme}

For the macroscopic homogenized system \eqref{TwoP}, we use the same type of scheme. The only difference is that the volume fraction of fluid $+$ does not satisfy $\alpha_+(1 - \alpha_+) = 0$ but
\[
D_t \alpha = \frac{\alpha (1 - \alpha)}{\alpha \mu_- + (1-\alpha) \mu_+}(p_+(\rho_+) - p_-(\rho_-) - (\mu_+ - \mu_-) \partial_x u).
\]
In the following we choose to discretize this equation in a forward Euler way (but a backward Euler scheme has also been tested and validated):
\begin{multline*}
\alpha_j^{n+1} = \alpha_j^n 
+ \Delta t \frac{\alpha_j^n(1 - \alpha_j^n)}{\alpha_j^n \mu_- + (1 - \alpha_j^n) \mu_+}\left(p_+(\rho_{+,j}^n) - p_-(\rho_{-,j}^n) - (\mu_+ - \mu_-) \frac{u_{j+1/2}^{n+1} - u_{j-1/2}^{n+1}}{x_{j+1/2}^{n+1} - x_{j-1/2}^{n+1}}\right).
\end{multline*}
All the other variables are approximated in a very standard and natural way.

\subsection{Experiments}

\begin{figure}[H]
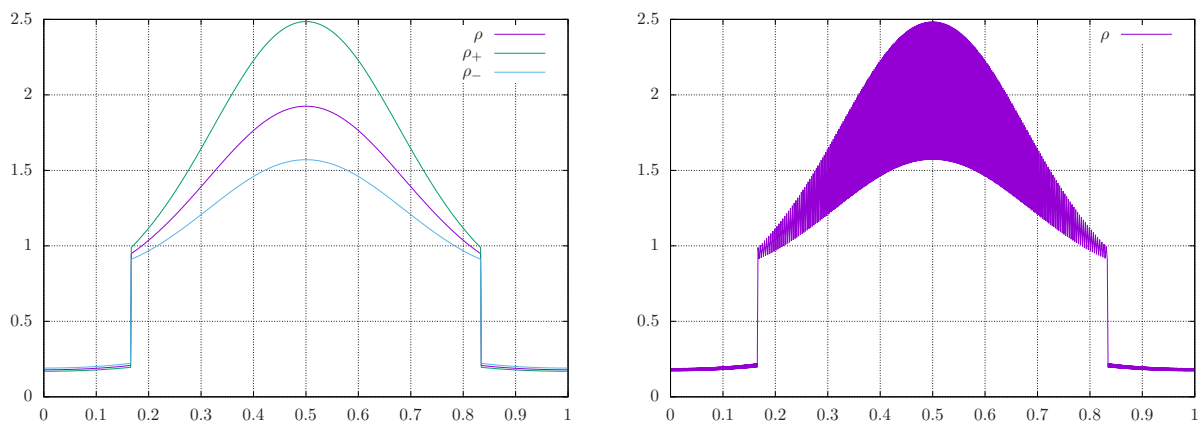

\begin{center}
\resizebox{\columnwidth}{!}
{\input{densities_macro.tex}
\input{density_meso.tex}}
\end{center}
\caption{Densities. On the left, the $3$ densities of the mixture, on the right, the density of the unmixed fluid. }\label{d}
\end{figure}

We propose two test-cases with $p_+(x) = x$ and $p_-(x) = x^2$. They are associated with a Cauchy datum of Riemann type:
\[
\left\{\begin{array}{l}
\alpha_0(x) = 1/2, \quad x \in \T_x, \\
\rho_+(x) = \rho_-(x) = \left\{\begin{array}{l}
1/8 \mbox{ if } x \in [0, 1/4) \cup [3/4, 1), \\
2 \mbox{ if } x \in [1/4, 3/4),
\end{array}\right. \\
u(x) = 0, \quad x \in \T_x, \\
\end{array}\right.
\]
and we propose to compare the numerical solutions obtained at time $t = 0.1$ with 1000 cells
\begin{itemize}
\item with the homogenized scheme of Section \ref{macroscheme},
\item and with the mesoscopic scheme of Section \ref{mesoscheme} by setting
\[
(\alpha_j^0, \rho_{+, j}^0, \rho_{-,j}^0) = \left\{\begin{array}{l}
(1,\rho_j^0, 0) \mbox{ if } j \mbox{ is even}, \\
(0,0, \rho_j^0) \mbox{ if } j \mbox{ is odd}, \\
\end{array}\right.
\]
and with a mesh with constant space step, which indeed corresponds in the weak limit to $\alpha = 1/2$. Note that the pressure is largely oscillating in this initial condition for the mesoscopic system.
\end{itemize}
In the first test, we take $\mu_+ = \mu_- = 0.1$ while in the second one we choose $\mu_+ = 0.1$ and $\mu_- = 0.02$.
Figures \ref{d} to \ref{vf} allow to compare the density, velocity, pressure and volume fraction. We observe a very good agreement between the mesoscopic and the macroscopic results. Note that for the mesoscopic  computation, we consider that there is only one density and one pressure, thus these quantities oscillate very fast (at the scale of the cell, which is the scale of the mixture). We observe, especially on the zoom of the density proposed by Figure \ref{dz}, that these oscillations occur between two functions that are very close to $\rho_+$ and $\rho_-$ computed by the macroscopic scheme. With the mesoscopic scheme, the volume fraction of fluid $+$ should oscillate between $0$ and $1$. In order to evaluate a volume fraction of $+$ in the limit mixture, what we here (Figure \ref{vf}) call $\alpha_j^n$ is computed by
\[
\alpha_j^n = \frac{c_j(x_{j+1/2}^n - x_{j-1/2}^n) + c_{j-1}(x_{j-1/2}^n - x_{j-3/2}^n)/2
+ c_{j+1}(x_{j+3/2}^n - x_{j+1/2}^n)/2}{x_{j+3/2}^n - x_{j-3/2}^n}
\]
(recall that $c_j$ is equal to $0$ or $1$ and does not depend on the time index).

\begin{figure}[H]
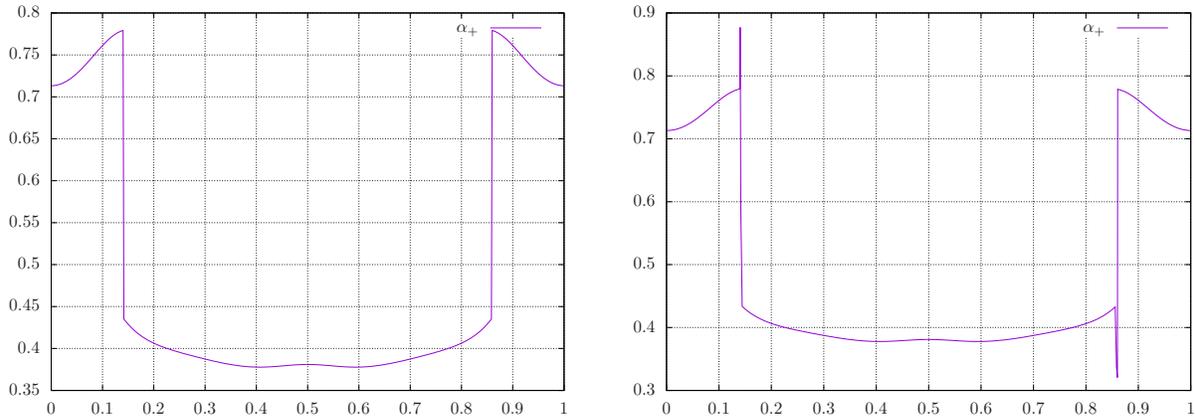

\begin{center}
\resizebox{\columnwidth}{!}
{\input{volume-fraction_macro_mudiff.tex}
\input{volume-fraction_meso_mudiff.tex}}
\caption{Volume fractions. On the left, the volume fraction $\alpha_+$ in the mixture, on the right, the estimate of the volume fraction in the unmixed fluid. }\label{vfm}
\end{center}
\end{figure}

The organization and the comments for the case with different viscosities, from Figure \ref{dm} to Figure \ref{vfm}, are the same.

\subsubsection{Case with equal viscosities}

\begin{figure}[H]
\begin{center}
\resizebox{0.7\width}{!}
{\input{densities_zoom.tex}}
\end{center}
\caption{Densities. Zoom of the preceding figures. }\label{dz}
\end{figure}

\begin{figure}[H]
\begin{center}
\resizebox{\columnwidth}{!}
{\input{velocity_macro.tex}
\input{velocity_meso.tex}}
\end{center}
\caption{Velocities. On the left, the velocity of the mixture, on the right, the velocity of the unmixed fluid. }\label{v}
\end{figure}

\begin{figure}[H]
\begin{center}
\resizebox{\columnwidth}{!}
{\input{pressures_macro.tex}
\input{pressure_meso.tex}}
\end{center}
\caption{Pressures. On the left, the $3$ pressures in the mixture, on the right, the pressure in the unmixed fluid. }\label{p}
\end{figure}

\begin{figure}[H]
\begin{center}
\resizebox{\columnwidth}{!}
{\input{volume-fraction_macro.tex}
\input{volume-fraction_meso.tex}}
\caption{Volume fractions. On the left, the volume fraction $\alpha_+$ in the mixture, on the right, the estimate of the volume fraction in the unmixed fluid. }\label{vf}
\end{center}
\end{figure}

\subsubsection{Case with different viscosities}

\begin{figure}[H]
\begin{center}
\resizebox{\columnwidth}{!}
{\input{densities_macro_mudiff.tex}
\input{density_meso_mudiff.tex}}
\end{center}
\caption{Densities. On the left, the $3$ densities of the mixture, on the right, the density of the unmixed fluid. }\label{dm}
\end{figure}

\begin{figure}[H]
\begin{center}
\resizebox{0.7\width}{!}
{\input{densities_zoom_mudiff.tex}}
\end{center}
\caption{Densities. Zoom of the preceding figures. }\label{dmz}
\end{figure}

\begin{figure}[H]
\begin{center}
\resizebox{\columnwidth}{!}
{\input{velocity_macro_mudiff.tex}
\input{velocity_meso_mudiff.tex}}
\end{center}
\caption{Velocities. On the left, the velocity of the mixture, on the right, the velocity of the unmixed fluid. }\label{vm}
\end{figure}

\begin{figure}[H]
\begin{center}
\resizebox{\columnwidth}{!}
{\input{pressures_macro_mudiff.tex}
\input{pressure_meso_mudiff.tex}}
\end{center}
\caption{Pressures. On the left, the $3$ pressures in the mixture, on the right, the pressure in the unmixed fluid. }\label{pm}
\end{figure}

\bigskip

\noindent {\bf Acknowledgments:}  The authors want to thank M. Hillairet and  L. Saint-Raymond 
for several discussions and comments. D. Bresch and C. Burtea are  partially 
supported bySingFlows project, grant ANR-18-CE40-002.

\end{document}